\documentstyle[leqno]{book}
\newfont{\cirilrm}{wncyr10 scaled 1000}
\newfont{\cirilbf}{wncyb10 scaled 1000}
\newfont{\cirilsf}{wncyss10 scaled 1000}
\newfont{\cirilit}{wncyi10 scaled 1000}
\newfont{\cirilsc}{wncysc10 scaled 1000}

\newcommand{\z}{\symbol{'31}}

\newcommand{\sh}{\symbol{'170}}

\newcommand{\ch}{\symbol{'161}}

\newcommand{\ja}{\symbol{'37}}

\newcommand{\mek}{\symbol{'176}}

\newcommand{\ij}{\symbol{'32}}

\newcommand{\ii}{\symbol{'171}}

\newcounter{supersection}[section]
\newtheorem{th}[supersection]{Theorem}
\newtheorem{lm}[supersection]{Lemma}
\newtheorem{re}[supersection]{Remark}
\newtheorem{co}[supersection]{Corollary}

\def\bibname{\textbf{REFERENCES}}
\def\thebibliography#1{\paragraph*{\uppercase{\bibname}}\list
{[\arabic{enumi}]}{\settowidth\labelwidth{[#1]}\leftmargin\labelwidth
\advance\leftmargin\labelsep\usecounter{enumi}}
\def\newblock{\hskip .11em plus .33em minus .07em}
\sloppy\clubpenalty4000\widowpenalty4000
\sfcode`\.=1000\relax}

\def\F11{{}_{1}\mbox{\rm F}{}_{1}}

\def\Dj{D{\hspace{-.75em}\raisebox{.3ex}{-}\hspace{.4em}}}
\def\stop{\mbox{\footnotesize {\vrule width 6pt height 6pt}}}

\begin{document}

\thispagestyle{plain}

\noindent $\;$ % {\small\sc ***. ***. ***. ***}

\noindent $\;$ % {\scriptsize ***. ***. (****), ??--??}

\vspace*{10.00 mm}

\centerline{\large \bf SOME CONSIDERATIONS IN CONNECTION}

\medskip

\centerline{\large \bf WITH KUREPA'S FUNCTION}
\footnotetext{2000 Mathematics Subject Classification: 30E20, 11J91.}
\footnotetext{Research partially supported by the MNTRS, Serbia \& Montenegro, Grant No. 1861.}

\vspace*{4.00 mm}

\centerline{\large \it Branko J. Male\v sevi\' c}

\vspace*{4.00 mm}

\begin{center}
\parbox{25.0cc}{\scriptsize \boldmath \bf
In this paper we consider the functional equation for factorial sum and
its particular solutions  (Kurepa's function $K(z)$ \cite{Kurepa_71} and
function $K_{1}(z)$). We determine an extension of domain of functions
$K(z)$ and $K_{1}(z)$ in the sense of Cauchy's principal value at point \cite{Slavic_70}.
In this paper we give an addendum to the proof of Slavi\' c's representation of
Kurepa's function $K(z)$ \cite{Slavic_73}. Also, we consider some representations of
functions $K(z)$ and $K_{1}(z)$ via incomplete gamma function and we consider differential
transcendency of previous functions too.}
\end{center}

\noindent
\section{\large \bf \boldmath \hspace*{-7.0 mm}
1. The functional equation for factorial sum  and its particular solutions}

The main object of consideration in this paper is the functional equation for
factorial sum:
\begin{equation}
\label{K_FE_1}
K(z) - K(z-1) = \Gamma(z),
\end{equation}
with respect to the function $K : D \longrightarrow C$ with domain
$D \subseteq C \backslash Z^{-}_{0}$, where $\Gamma$ is the gamma
function, $C$ is the set of complex numbers and $Z^{-}_{0}$ is the
set of negative integer numbers. A solution of functional equation
(\ref{K_FE_1}) over the set of natural numbers ($D = N$)
is the function of left factorial $!n$. {\sc \Dj uro Kurepa}
introduced this function, in the paper \cite{Kurepa_71},
as sum of factorials $!n = 0! + 1! + 2! + \ldots + (n-1)!$.
Let us use the standard notation:
\begin{equation}
\label{K_SUM_1}
K(n) = \displaystyle\sum\limits_{i=0}^{n-1}{i!}.
\end{equation}
Sum (\ref{K_SUM_1}) corresponds to the sequence $A003422$ in
\cite{Sloane_03}. We call the functional equation (\ref{K_FE_1})
{\em the functional equation for factorial sum}. In consideration
which follows we consider two particular solutions of the
functional equation (\ref{K_FE_1}).

\bigskip
\noindent {\bf \boldmath 1.1. The function $K(z)$.} An analytical
extension of the function (\ref{K_SUM_1}) over the set of complex
numbers is determined by integral~\cite{Kurepa_73}:
\begin{equation}
\label{K_INT_1}
K(z)
=
\displaystyle\int\limits_{0}^{\infty}{
e^{-t} \displaystyle\frac{t^{z}-1}{t-1} \: dt},
\end{equation}
which converges for $\mbox{Re} \: z > 0$. For the function $K(z)$
we use the term Kurepa's function and it is a solution of the
functional equation (\ref{K_FE_1}). Let us observe that since
\mbox{$K(z-1) = K(z) - \Gamma(z)$}, it is possible to make
analytical continuation of Kurepa's function $K(z)$ for
\mbox{$\mbox{Re} \, z \leq 0$}. In that way, the Kurepa's function
$K(z)$ is a meromorphic function with simple poles at $z=-1$ and
$z=-n$ $(n \!\geq\! 3)$. At point $z = -2$ Kurepa's function has a
removable singularity and \mbox{$K(-2) \stackrel{\mbox{\tiny def}}{=}
\lim_{z \rightarrow -2}{K(z)} = 1$}. Kurepa's function has an essential
singularity at point $z = \infty$. Kurepa's function has the following residues:
\begin{equation}
\begin{array}{l}
\mathop{\mbox{\rm res}}\limits_{z = -1}{K(z)} = -1, \\
\mathop{\mbox{\rm res}}\limits_{z = -n}{K(z)} =
\displaystyle\sum\limits_{k=2}^{n-1}{\displaystyle\frac{(-1)^{k-1}}{k!}}
\quad (n\!\geq\!3).
\end{array}
\end{equation}
Previous results for Kurepa's function are presented according to \cite{Kurepa_73}
and \cite{Slavic_73}.

\bigskip
\noindent
{\bf \boldmath 1.2. The function $K_{1}(z)$.}
The functional equation (\ref{K_FE_1}), besides Kurepa's function $K(z)$,
has another solution which is given by the following statement.
\begin{th}
\label{K_lema_4}
Let $D = C \backslash Z$. Then, series$:$
\begin{equation}
\label{Def_K1}
K_{1}(z) = \displaystyle\sum\limits_{n=0}^{\infty}{\Gamma(z-n)}
\end{equation}
absolutely converges and it is a solution of the functional equation
{\rm (\ref{K_FE_1})}~over~$D$.
\end{th}

\noindent {\bf Proof.} Set $D$ has the following decomposition in
two disjoint sets:
\begin{equation}
D_{1} = \{ z \!\in\! D \; | \; \mbox{\rm Re}\,z \not\in Z \}
\end{equation}
and
\begin{equation}
D_{2}
=
\{ z \!\in\! D \; | \; \mbox{\rm Re}\,z \in Z \; \wedge \; \mbox{\rm Im}\,z \neq 0 \}.
\end{equation}

\noindent We prove the statement by  discussing the following two
cases.

\smallskip \noindent
{\boldmath $1^{0}.$} Let  $z \!=\! x\!+\!iy \!\in\! D_{1}$. Let us
denote $m \!=\! [x]$. For each $n\!\in\!N_{0}$, such that
\mbox{$n\!\geq\!m\!+\!2$}, it is true that:
\begin{equation}
|\Gamma(z\!-\!n)|
\leq
|\Gamma(x\!-\!n)|
=
\displaystyle\frac{\pi}{\, {\big |} \sin {\big (} \pi x {\big )}
\, \Gamma {\big (} n\!+\!1\!-\!x {\big )} {\big |}}
<
\displaystyle\frac{\pi}{\, | \sin ( \pi x ) |}
\cdot
\displaystyle\frac{1}{\,\Gamma(n\!-\!m)}.
\end{equation}
Thus, for $m_{1} = \max \{m\!+\!2, 0\}$, the following is true:
\begin{equation}
\displaystyle\sum\limits_{n=m_{1}}^{\infty}{|\Gamma(z\!-\!n)|}
<
\displaystyle\frac{\pi}{\, | \sin ( \pi x ) |}
\cdot
\displaystyle\sum\limits_{n=m_{1}}^{\infty}{\displaystyle\frac{1}{\,\Gamma(n\!-\!m)}},
\end{equation}
which is sufficient to conclude that the statement is true over $D_{1}$.

\smallskip \noindent
{\boldmath $2^{0}.$} Let  $z \!=\! x+iy \!\in\! D_{2}$. Let us
denote by $m\!=\!\mbox{\rm Re}\,z$. For each $n\!\in\!N_{0}$, such
that \mbox{$n\!\geq\!m\!+\!1$}, it is true that:
\begin{equation}
\begin{array}{rcl}
|\Gamma(z\!-\!n)|
&\!\!=\!\!&
|\Gamma(iy - (n\!-\!m))|
=
\displaystyle\frac{{\big |} \Gamma(iy) {\big |}}{\,{\big |}
\displaystyle\prod\limits_{r=1}^{n-m}{\!(-r\!+\!iy)}\,{\big |}\,}
=
\displaystyle\frac{\sqrt{\displaystyle\frac{\pi}{\,y \, \mbox{\rm sh}(\pi y)}}}{\,
\displaystyle\prod\limits_{r=1}^{n-m}{|\!-\!r\!+\!iy \,|}}                     \\[7.5 ex]
&\!\!\!=\!\!\!&
\sqrt{\displaystyle\frac{\pi}{\,y \, \mbox{\rm sh}(\pi y)}} \cdot
\displaystyle\frac{1}{\,\displaystyle\prod\limits_{r=1}^{n-m}{\!\sqrt{r^2\!+\!y^2}}\,}
\; < \;
\sqrt{\displaystyle\frac{\pi}{\,y \, \mbox{\rm sh}(\pi y)}} \cdot
\displaystyle\frac{1}{\,\Gamma(n\!-\!m\!+\!1)} \, .
\end{array}
\end{equation}
Thus, for $m_{2} = \max \{m\!+\!1, 0\}$, the following is true:
\begin{equation}
\displaystyle\sum\limits_{n=m_{2}}^{\infty}{|\Gamma(z\!-\!n)|}
<
\sqrt{\displaystyle\frac{\pi}{\,y \, \mbox{\rm sh}(\pi y)}} \cdot
\displaystyle\sum\limits_{n=m_{2}}^{\infty}{\displaystyle\frac{1}{\,\Gamma(n\!-\!m\!+\!1)}}\,,
\end{equation}
which is sufficient to conclude that the statement is true over
$D_{2}\!=\!D \backslash D_{1}$. Let us note that it is easy to
prove that the function $K_{1}(z)$ is a solution of the functional
equation (\ref{K_FE_1}).~\stop
\begin{re}
Function $K_{1}(z)$, defined by {\rm (\ref{Def_K1})} over $C$,
has poles at integer points $z\!=\!m \in Z$.
\end{re}

\section{\large \bf \boldmath \hspace*{-7.0 mm}
2. Extending the domain of functions $K(z)$ and $K_{1}(z)$ in the
\hspace*{0.0 mm}sense of Cauchy's  principal value}

Let us observe a possibility of extending the domain of the functions $K(z)$ and $K_{1}(z)$,
in the sense of Cauchy's principal value, over the set of complex numbers. Namely, for
a meromorphic function $f(z)$, on the basis of Cauchy's integral formula, we define
{\em the principal value at point $a$} as follows~\cite{Slavic_70}:
\begin{equation}
\label{GAMMA_PV_a}
\mathop{\mbox{\rm p.v.}}\limits_{z = a}{f(z)}
=
\lim\limits_{\rho \rightarrow 0_{+}}{\displaystyle\frac{1}{2 \pi i}
\!\!\!\!\displaystyle\oint\limits_{|z-a|=\rho}{\!\!\!\!\displaystyle\frac{f(z)}{z-a}\,dz}}.
\end{equation}
It is easily proved that the principal value at pole $z = a$ exists as a finite
complex number. For two meromorphic functions $f_1(z)$ and $f_2(z)$ additivity
is true~\cite{Slavic_70}:
\begin{equation}
\mathop{\mbox{\rm p.v.}}\limits_{z = a}{{\Big (}f_{1}(z)+f_{2}(z){\Big )}}
=
\mathop{\mbox{\rm p.v.}}\limits_{z = a}{f_{1}(z)}
+
\mathop{\mbox{\rm p.v.}}\limits_{z = a}{f_{2}(z)}.
\end{equation}
The main formula in this section is given by the following statement:
\begin{lm}
\label{TH_Gamma_PV_f}
For function $f(z)$ with simple pole at point $z = a$ the following is true$:$
\begin{equation}
\label{PV_f_1}
\mathop{\mbox{\rm p.v.}}\limits_{z = a}{f(z)}
=
\lim\limits_{\varepsilon \rightarrow 0}{
\displaystyle\frac{f(a-\varepsilon)+f(a+\varepsilon)}{2}}.
\end{equation}
\end{lm}

\noindent
{\bf Proof.} Let $z=a$ be a simple pole of the function $f(z)$. Then, there
exist $\rho > 0$ such that function $f(z)$ has the following representation:
$f(z) = \displaystyle \frac{g(z)}{z-a}$, for some regular function $g(z)$ over
$|z-a|<\rho$. Hence, using Cauchy's integral formula:
\begin{equation}
\mathop{\mbox{\rm p.v.}}\limits_{z = a}{f(z)}
=
\lim\limits_{\rho \rightarrow 0}{\displaystyle\frac{1}{2 \pi i}
\!\!\!\!\displaystyle\oint\limits_{|z-a|=\rho}{\!\!\!\!\displaystyle\frac{f(z)}{z-a}\,dz}}
=
\lim\limits_{\rho \rightarrow 0}{\displaystyle\frac{1}{2 \pi i}
\!\!\!\!\displaystyle\oint\limits_{|z-a|=\rho}{\!\!\!\!\displaystyle\frac{g(z)}{(z-a)^2}\,dz}}
=
g'(a).
\end{equation}
On the other hand:
\begin{equation}
\lim\limits_{\varepsilon \rightarrow 0}{
\displaystyle\frac{f(a-\varepsilon)+f(a+\varepsilon)}{2}}
=
\lim\limits_{\varepsilon \rightarrow 0}{
\displaystyle\frac{g(a+\varepsilon)-g(a-\varepsilon)}{2 \varepsilon}}
=
g'(a).
\end{equation}
Previous two equalities are sufficient that we conclude that (\ref{PV_f_1}) is true.~\stop
\begin{co}
\label{TH_Gamma_PV}
For gamma function $\Gamma(z)$ it is true$:$
\begin{equation}
\label{PV_1}
\mathop{\mbox{\rm p.v.}}\limits_{z = -n}{\!\!\Gamma(z)}
\!=\!
\lim\limits_{\varepsilon \rightarrow 0}{\!\!
\displaystyle\frac{\Gamma(-n-\varepsilon)+\Gamma(-n+\varepsilon)}{2}}
\!=\!
(-1)^{n}\displaystyle\frac{\Gamma^{'}(n+1)}{\Gamma(n+1)^2}
\quad (n\!\in\!N_{0}).
\end{equation}
\end{co}

\noindent
{\bf Proof.} On the basis of equality $\Gamma(-z) \!=\!
\displaystyle\frac{-\pi}{\Gamma(z+1) \sin \pi z}$, $(z \!\not\in\! Z)$,
for real $\varepsilon\!>\!0$, is true:
\begin{equation}
\displaystyle\frac{\Gamma(-n - \varepsilon) + \Gamma(-n + \varepsilon)}{2}
=
(-1)^{n-1}
\displaystyle\frac{\pi \varepsilon}{\sin \pi \varepsilon}
\cdot
\displaystyle\frac{
\mbox{\small $\displaystyle\frac{1}{\Gamma(n\!+\!1\!+\!\varepsilon)}$}
-
\mbox{\small $\displaystyle\frac{1}{\Gamma(n\!+\!1\!-\!\varepsilon)}$}}{2 \varepsilon}.
\end{equation}
Hence, using real limit $\varepsilon \rightarrow 0_{+}$, we obtain result (\ref{PV_1})
from \cite{Slavic_70} in an easier way.~\stop
\begin{re}
\label{Gamma_Kolicnik_Je_Opadajuci_Niz}
For $n \in N_{0}$ it is true~{\rm \cite{Slavic_70}}$:$

\vspace*{-2.5 mm}

\begin{equation}
\displaystyle\frac{\Gamma^{'}(n+1)}{\Gamma(n+1)^2}
=
\displaystyle\frac{\mbox{\small $-\gamma + 1 + \displaystyle\frac{1}{2} + \ldots +
\displaystyle\frac{1}{n}$}}{n!},
\end{equation}
where $\gamma$ is Euler's constant.
\end{re}
Extension of the domain of the functions $K(z)$ and $K_{1}(z)$, in the sense of
Cauchy's principal value, is given by the following two theorems.
\begin{th}
For Kurepa's function $K(z)$ it is true$:$
\begin{equation}
\label{K_PV}
\mathop{\mbox{\rm p.v.}}\limits_{z = -n}{\!\!K(z)}
=
-\displaystyle\sum\limits_{i=0}^{n-1}{
\mathop{\mbox{\rm p.v.}}\limits_{z=-i}{\Gamma(z)}}
=
\displaystyle\sum\limits_{i=0}^{n-1}{(-1)^{i+1}
\displaystyle\frac{\Gamma^{'}(i+1)}{\Gamma(i+1)^2}}
\quad (n \!\in\! N).
\end{equation}
\end{th}

\noindent
{\bf Proof.} If equality:
\begin{equation}
K(z) = K(z+n) - {\big (} \Gamma(z+1) + \ldots + \Gamma(z+n) {\big )}
\end{equation}
we consider at the point $z = -n$ in the sense of Cauchy's principal value,
on the basis of (\ref{PV_1}), the equality (\ref{K_PV}) follows. Let us remark
that $\mathop{\mbox{\rm p.v.}}\limits_{z = -2}{\!\!K(z)} = K(-2) = 1$.~\stop
\begin{lm}
Let's define $L_{1}
=
-\displaystyle\sum\limits_{n=0}^{\infty}{
\!\,\mathop{\mbox{\rm p.v.}}\limits_{z = -n}{\!\,\Gamma(z)} }
$, then:
\begin{equation}
\label{Const_L_1}
L_{1}
=
\displaystyle\sum\limits_{n=0}^{\infty}{
(-1)^{n+1} \displaystyle\frac{\Gamma^{'}(n\!+\!1)}{\Gamma(n\!+\!1)^2} }
\approx
0.697 \, 174 \, 883 \, .
\end{equation}
\end{lm}

\noindent
{\bf Proof.} The previous series converges on the basis of the remark
\ref{Gamma_Kolicnik_Je_Opadajuci_Niz}.
\begin{th}
\label{K_aux}
For function $K_{1}(z)$ it is true$:$
\begin{equation}
\mathop{\mbox{\rm p.v.}}\limits_{z = n}{K_{1}(z)}
=
\mathop{\mbox{\rm p.v.}}\limits_{z = n}{K(z)} - L_{1}
\quad (n \!\in\! Z).
\end{equation}
\end{th}

\noindent
{\bf Proof.} For $n \geq 0$ it is true:
\begin{equation}
\mathop{\mbox{\rm p.v.}}\limits_{z = n}{K_{1}(z)}
\!=\!
\displaystyle\sum\limits_{i=0}^{\infty}{
\!\mathop{\mbox{\rm p.v.}}\limits_{z = n-i}{\!\Gamma(z)}}
\!=\!
\displaystyle\sum\limits_{i=0}^{\infty}{
\!\mathop{\mbox{\rm p.v.}}\limits_{z = -i}{\!\Gamma(z)}}
+
\displaystyle\sum\limits_{i=1}^{n}{\Gamma(i)}
\!=\!
K(n)
-
L_{1}.
\end{equation}
For $n < 0$ it is true:
\begin{equation}
\begin{array}{rcl}
\mathop{\mbox{\rm p.v.}}\limits_{z = n}{K_{1}(z)}
&\!\!\!=\!\!\!&
\displaystyle\sum\limits_{i=0}^{\infty}{
\!\mathop{\mbox{\rm p.v.}}\limits_{z = n-i}{\!\Gamma(z)}}
\!=\!
\displaystyle\sum\limits_{i=0}^{\infty}{
\!\mathop{\mbox{\rm p.v.}}\limits_{z = -i}{\!\Gamma(z)}}
\; - \!
\displaystyle\sum\limits_{i=0}^{(-n)-1}{
\!\mathop{\mbox{\rm p.v.}}\limits_{z = -i}{\!\Gamma(z)}}                       \\[2.0 ex]
&\!\!\!=\!\!\!&
{\bigg (} - \!\displaystyle\sum\limits_{i=0}^{(-n)-1}{
\mathop{\mbox{\rm p.v.}}\limits_{z = -i}{\Gamma(z)}} {\bigg )}
-
L_{1}
=
\mathop{\mbox{\rm p.v.}}\limits_{z = n}{K(z)}
-
L_{1}.\;\;\stop
\end{array}
\end{equation}

\noindent
\section{\large \bf \boldmath \hspace*{-7.0 mm}
3. Slavi\' c's formula for Kurepa's function}

\medskip
\noindent In this section we give an addendum to the proof of
Slavi\' c's representation of Kurepa's function $K(z)$
\cite{Slavic_73}, by the following two statements:
\begin{lm}
\label{K_lema_12}
Function$:$
\begin{equation}
\label{Maple_K_Lema_12}
F(z)
=
\displaystyle\sum\limits_{n=1}^{\infty}{
\!{\bigg (}\displaystyle\sum\limits_{k=1}^{\infty}{
\displaystyle\frac{(-1)^{n+k-1}(n+k+1)}{(n+k)!} \, z^k {\bigg )}}},
\end{equation}
is entire, whereas following is true$:$
\begin{equation}
\label{Maple_K_Lema_21}
F(z)
=
\displaystyle\sum\limits_{k=1}^{\infty}{
\!{\bigg (}\displaystyle\sum\limits_{n=1}^{\infty}{
\displaystyle\frac{(-1)^{n+k-1}(n+k+1)}{(n+k)!} \, z^k {\bigg )}}}
=
e^{-z}-1.
\end{equation}
\end{lm}

\noindent {\bf Proof.} For $z = 0$ the equality (\ref{Maple_K_Lema_21}) is true. Let us
introduce a sequence of functions:
\begin{equation}
\label{Def_f_n}
f_{n}(z)
=
\displaystyle\sum\limits_{k=1}^{\infty}{
\displaystyle\frac{(-1)^{n+k-1}(n\!+\!k\!+\!1)}{(n\!+\!k)!} \, z^k},
\end{equation}
for $z \in C$ $(n \in N)$. Previous series converge over C
because, for $z \neq 0$, it is true that:
\begin{equation}
\label{Maple_f_n}
f_{n}(z)
=
\displaystyle\sum\limits_{j=0}^{n}{
(-1)^{j} {\bigg (}\displaystyle\frac{j}{j!}
\!+\!
\displaystyle\frac{1}{j!} {\bigg )} \, z^{j-n}
\!+\!
e^{-z}(z^{-n+1}\!-\!z^{-n})}.
\end{equation}

\vspace*{-2.5 mm}

\noindent Let us mention that the previous equality is easily
checked by the following substitution \mbox{$e^{-z} =
\sum_{k=0}^{\infty}{\frac{(-z)^k}{k!}}$} at the right side of
equality of formula (\ref{Maple_f_n}). Let  $\rho>0$ be fixed.
Over the set \mbox{$D\!=\!\{ z\!\in\!C \, | \,
0\!<\!|z|\!<\!\rho\}$} let us an auxiliary function \mbox{$g(z) =
(z-1)e^{-z} : D \longrightarrow C$}. If we denote by $R_{n}(.)$
the remainder of $n$-th order of MacLaurin's expansion, then for
$z \in D$ the following representation is true:
\begin{equation}
\label{MacLaurin_f_n}
f_{n}(z)
=
\displaystyle\frac{R_{n}(g(z))}{z^n}.
\end{equation}
Over $E=(0,\rho)$ let us form an auxiliary function \mbox{$h(t) =
(t+1)e^{t} : E \longrightarrow R^{+}$}. Then, for $|z| < \rho$ it
is true:
\begin{equation}
\label{Ineq_f_n}
|f_{n}(z)|
\leq
\displaystyle
\frac{e^{\rho}(n\!+\!2\!+\!\rho)}{n!}\rho.
\end{equation}
Indeed, previous inequality is true for $z=0$. For $z \in D$ and $t = |z| \in E$
exists $c \in (0,t)$ such that:
\begin{equation}
\quad |f_{n}(z)|
=
{\Bigg |}\displaystyle\frac{R_{n}(g(z))}{z^n}{\Bigg |}
\leq
\displaystyle\frac{R_{n}(h(t))}{t^n}
=
\displaystyle
\frac{h^{(n+1)}(c)}{(n+1)!}\,t
\leq
\displaystyle
\frac{e^{\rho}(n\!+\!2\!+\!\rho)}{n!}\rho.
\end{equation}
For function:
\begin{equation}
F(z)
=
\displaystyle
\sum\limits_{n=1}^{\infty}{f_{n}(z)},
\end{equation}
it is possible, for $|z|\!<\!\rho$, to apply Weierstrass's double series
theorem~\cite{Knopp_96} (page~83). Indeed, on the basis of (\ref{Def_f_n}),
the functions $f_{n}(z)$ are regular for $|z|\!<\!\rho$. On the basis of
(\ref{Ineq_f_n}), the series $\sum_{n=1}^{\infty}{f_{n}(z)}$ is uniformly
convergent for $|z|\!\leq\!r\!<\!\rho$, for every $r\!<\!\rho$. Then on
the basis of the Weierstrass's theorem, for $|z|\!<\!\rho$, the following~is~true:
\begin{equation}
F(z)
=
\displaystyle\sum\limits_{k=1}^{\infty}{
\!{\bigg (}\displaystyle\sum\limits_{n=1}^{\infty}{
\displaystyle\frac{(-1)^{n+k-1}(n+k+1)}{(n+k)!} \, z^k {\bigg )}}}
=
e^{-z}-1,
\end{equation}
because:
\begin{equation}
\label{Razmena}
\label{K_sum} \displaystyle\sum\limits_{n=1}^{\infty}{
\displaystyle\frac{(-1)^{n+k-1}(n\!+\!k\!+\!1)}{(n\!+\!k)!}} =
\displaystyle\frac{(-1)^{k}}{k!}.
\end{equation}
The previous series is a telescoping series. Let us note that $\rho>0$ can be arbitrarily
large positive number. Hence, the equality (\ref{Maple_K_Lema_21}) is true for all $z \in C$;
i.e. the function $F(z)$ is entire. \stop
\begin{re}
For $|z| > 1$ the following equalities are true :
\begin{equation}
\label{Maple_K_Sum_1}
\displaystyle\sum\limits_{n=1}^{\infty}{
e^{-z}(z^{-n+1}\!-\!z^{-n})} = e^{-z}
\end{equation}
and
\begin{equation}
\label{Maple_K_-_1}
\displaystyle\sum\limits_{n=1}^{\infty}{\!{\bigg (}
\displaystyle\sum\limits_{j=0}^{n}{(-1)^j\displaystyle\frac{j}{j!} \, z^{j-n}}{\bigg )}}
=
-1
-
\displaystyle\sum\limits_{n=1}^{\infty}{\!{\bigg (}
\displaystyle\sum\limits_{j=0}^{n}{(-1)^{j}\displaystyle\frac{1}{j!} \,
z^{j-n}}{\bigg )}}.
\end{equation}
Hence for $|z| > 1$, on the basis of {\rm (\ref{Maple_f_n})}, it also follows
that {\rm (\ref{Maple_K_Lema_21})} is true.
\end{re}
\begin{lm}
\label{K_lema_3}
For $z \in C$ it is true$:$
\begin{equation}
\label{K_sum_sum}
(z\!-\!1)
\displaystyle\sum\limits_{n=1}^{\infty}{
\displaystyle\sum\limits_{k=0}^{\infty}{
\displaystyle\frac{(-1)^{k+n}}{(k\!+\!n)!} \, z^{k}}}
=
e^{-z} \!- e^{-1}.
\end{equation}
\end{lm}

\noindent {\bf Proof.} On the basis of the lemma \ref{K_lema_12}
it is true that:
\begin{equation}
\mbox{\small $
\begin{array}{rcl}
(z\!-\!1)
\displaystyle\sum\limits_{n=1}^{\infty}{
\displaystyle\sum\limits_{k=0}^{\infty}{
\displaystyle\frac{(-1)^{k+n}}{(k\!+\!n)!} \, z^{k}}}
&\!\!=\!\!&
\displaystyle\sum\limits_{n=1}^{\infty}{
\displaystyle\sum\limits_{k=0}^{\infty}{
\displaystyle\frac{(-1)^{k+n}}{(k\!+\!n)!} \, (z^{k+1}\!-z^{k})}}              \\[2.5 ex]
&\!\!=\!\!&
\displaystyle\sum\limits_{n=1}^{\infty}{
\!{\bigg (}\!\!-\displaystyle\frac{(-1)^n}{n!}
+\!\displaystyle\sum\limits_{k=1}^{\infty}{
{\Big (}\displaystyle\frac{(-1)^{k+n-1}}{(k\!+\!n\!-\!1)!}
\!-\!\displaystyle\frac{(-1)^{k+n}}{(k\!+\!n)!}
{\Big )}z^{k}}{\bigg )}}                                                       \\[2.5 ex]
&\!\!\mathop{=}\limits_{(\ref{Maple_K_Lema_21})}\!\!&
\displaystyle\sum\limits_{k=1}^{\infty}{
{\bigg (}\displaystyle\sum\limits_{n=1}^{\infty}{
\displaystyle\frac{(-1)^{n+k-1}(n\!+\!k\!+\!1)}{(n\!+\!k)!}}{\bigg )}\,z^{k}}
-\displaystyle\sum\limits_{n=1}^{\infty}{
\displaystyle\frac{(-1)^n}{n!}}                                                \\[5.0 ex]
&\!\!\mathop{=}\limits_{(\ref{Maple_K_Lema_21})}\!\!&
e^{-z} - e^{-1}.\;\stop
\end{array}$}
\end{equation}

\smallskip \noindent
On the basis of the previous lemmas we give an addendum to the
proof of Slavi\' c's formula for Kurepa's
function~\cite{Slavic_73}:
\begin{th}
For Kurepa's function $K(z)$ the following representation is true$:$
\begin{equation}
\label{K_sum_Slavic}
K(z)
=
\displaystyle\frac{\mbox{\rm Ei}(1)}{e}
-
\displaystyle\frac{\pi}{e} \, \mbox{\rm ctg} \, \pi z
+
\displaystyle\sum\limits_{n=0}^{\infty}{\Gamma(z-n)},
\end{equation}
where the values in the previous formula, in integer points $z$,
are determined in the sense of Cauchy's principal value.
\end{th}

\noindent
{\bf Proof.} Some parts of this proof are reproduced according to \cite{Slavic_73}.
For $-(n\!+\!1)\!<\!\mbox{\rm Re}\:z\!<\!-n$ and $n=0,1,2,\ldots$ the following
formula is true \cite{Bateman_65}:
\begin{equation}
\Gamma(z)
=\!\!
\displaystyle\int\limits_{0}^{+\infty}\!{{\bigg (}\!
e^{-t} - \displaystyle\sum\limits_{m=0}^{n}{
\displaystyle\frac{(-t)^{m}}{m!}}\!{\bigg )}t^{z-1}\:dt}.
\end{equation}
Hence, for $0 < \mbox{\rm Re}\:z < 1$ and $n=1,2,\ldots$ the following formula is true:

\vspace*{-2.5 mm}

\begin{equation}
\Gamma(z-n)
=\!\!
\displaystyle\int\limits_{0}^{+\infty}\!{{\bigg (}\!
e^{-t} - \displaystyle\sum\limits_{m=0}^{n-1}{
\displaystyle\frac{(-t)^{m}}{m!}}\!{\bigg )}t^{z-n-1}\:dt}.
\end{equation}
Further we observe the following difference:
\begin{equation}
\begin{array}{rcl}
K(z) - \displaystyle\sum\limits_{n=0}^{+\infty}{\Gamma(z\!-\!n)}
&\!\!=\!\!&
\!\displaystyle\int\limits_{0}^{+\infty}{\!
e^{-t}\displaystyle\frac{t^z\!-\!1}{t\!-\!1}\:dt}
-\!
\displaystyle\int\limits_{0}^{+\infty}{\!
e^{-t}t^{z\!-\!1}\:dt}                                                         \\[2.0 ex]
&\!\!-\!\!&
\!\displaystyle\sum\limits_{n=1}^{\infty}{\!
\displaystyle\int\limits_{0}^{+\infty}{\!\!{\bigg (}\!
e^{-t}\!-\!\displaystyle\sum\limits_{m=0}^{n-1}{\!
\displaystyle\frac{(-t)^{m}}{m!}}\!{\bigg )}t^{z\!-\!n\!-\!1}\:dt}}.
\end{array}
\end{equation}
In this part of the proof we give an addendum using lemma \ref{K_lema_3}.
The following derivation is true:
\begin{equation}
\label{K_central}
\mbox{\small $\begin{array}{rcl}
K(z) - \displaystyle\sum\limits_{n=0}^{+\infty}{\Gamma(z-n)}
\!\!&\!\!=\!\!&\!\!
\!\displaystyle\int\limits_{0}^{\infty}{\!
e^{-t} \displaystyle\frac{t^{z-1}\!-\!1}{t\!-\!1} \: dt}
-\!
\displaystyle\int\limits_{0}^{+\infty}{\!
\displaystyle\sum\limits_{n=1}^{\infty}{\!\!
{\bigg (}\!e^{-t}\!-\!\displaystyle\sum\limits_{m=0}^{n-1}{\!
\displaystyle\frac{(-t)^{m}}{m!}}\!{\bigg )}t^{z\!-\!n\!-\!1}\:dt}}            \\[2.0 ex]
\!\!&\!\!=\!\!&\!\!
\!\displaystyle\int\limits_{0}^{\infty}{
{\bigg (}
\displaystyle\frac{e^{-t}}{1\!-\!t}
+
e^{-t}\displaystyle\frac{t^{z-1}}{t\!-\!1}
-
\displaystyle\sum\limits_{n=1}^{\infty}{
\displaystyle\sum\limits_{m=n}^{\infty}{
\displaystyle\frac{(-t)^{m}}{m!}}\,t^{z\!-\!n\!-\!1}}{\bigg )}dt}              \\[2.0 ex]
\!\!&\!\!=\!\!&\!\!
\!\displaystyle\int\limits_{0}^{\infty}{
{\bigg (}
\displaystyle\frac{e^{-t}}{1\!-\!t}
+
\displaystyle\frac{t^{z-1}}{t\!-\!1}
\!{\bigg (}\!e^{-t}-\!(t\!-\!1)\!
\displaystyle\sum\limits_{n=1}^{\infty}{
\displaystyle\sum\limits_{m=n}^{\infty}{
\displaystyle\frac{(-1)^{m}}{m!}\,t^{m-n}}}\!{\bigg )}\!\!{\bigg )}dt}         \\[2.0 ex]
\!\!&\!\!=\!\!&\!\!
\!\displaystyle\int\limits_{0}^{\infty}{\!
{\bigg (}
\displaystyle\frac{e^{-t}}{1\!-\!t}
+
\displaystyle\frac{t^{z-1}}{t\!-\!1}
\!{\bigg (}\!e^{-t}-\!(t\!-\!1)\!
\displaystyle\sum\limits_{n=1}^{\infty}{
\displaystyle\sum\limits_{k=0}^{\infty}{
\displaystyle\frac{(-1)^{k+n}}{(k\!+\!n)!}\,t^{k}}}\!{\bigg )}\!\!{\bigg )}dt} \\[2.0 ex]
&\!\!\mathop{=}\limits_{(\ref{K_sum_sum})}\!\!&
\!\displaystyle\int\limits_{0}^{\infty}{\!
{\bigg (}
\displaystyle\frac{e^{-t}}{1\!-\!t}
+
\displaystyle\frac{1}{e}
\displaystyle\frac{t^{z-1}}{t\!-\!1} {\bigg )} dt},
\end{array}$}
\end{equation}

\vspace*{-2.5 mm}

\noindent
for $0 < \mbox{\rm Re}\:z < 1$. Integral at the right side of equality (\ref{K_central}),
which converges in ordinary sense, will be substituted by two integrals which converge
in the sense of the Cauchy's principal value. Namely, using the function of exponential
integral~\cite{Bateman_65}:
\begin{equation}
\label{Ei}
\mbox{\rm Ei}(x) =
\mbox{\rm p.v.}\!\displaystyle\int\limits_{-\infty}^{x}{\!\displaystyle\frac{e^t}{t} \, dt},
\end{equation}
and using the formulas 3.352-6 and 3.311-8 from \cite{Rizik_71}:
\begin{equation}
\mbox{\rm p.v.} \, \displaystyle\int\limits_{0}^{\infty}{
\displaystyle\frac{e^{-t}}{1-t} \, dt} = \displaystyle\frac{\mbox{Ei}(1)}{e}
\end{equation}
and
\begin{equation}
\mbox{\rm p.v.} \, \displaystyle\int\limits_{0}^{\infty}{
\displaystyle\frac{t^{z-1}}{1-t} \, dt} = \pi \, \mbox{\rm ctg} \, \pi z
\end{equation}

\break

\noindent
we  conclude that formula (\ref{K_sum_Slavic}) is true for $0 < \mbox{\rm Re}\:z < 1$.
According to Riemann's theorem we conclude that Slavi\' c's formula (\ref{K_sum_Slavic}) is
true for each complex $z$. Namely, formula (\ref{K_sum_Slavic}), in integer points $z$, is
true in the sense of Cauchy's principal value on the basis of the theorem \ref{K_aux}.$\;\stop$
\begin{co}
Function $K_{1}(z)$ is a meromorphic function with simple poles in
integer points $z=m$ $(m\!\in\!Z)$ and with residue values$:$
\begin{equation}
\mathop{\mbox{\rm res}}\limits_{z = m}{K_{1}(z)}
=
\displaystyle\frac{1}{e} + \mathop{\mbox{\rm res}}\limits_{z = m}{K(z)} \quad (m\!\in\!Z).
\end{equation}
At the point $z = \infty$ function $K_{1}(z)$ has an essential
singularity.
\end{co}

\noindent
\section{\large \bf \boldmath \hspace*{-7.0 mm}
4. Some representations of functions $K(z)$ and $K_{1}(z)$ via
\hspace*{3.5 mm}incomplete gamma function}

In this section we give some representations of functions $K(z)$ and $K_{1}(z)$ via
gamma and incomplete gamma functions, where the later ones % (incomplete gamma functions)
are defined by integrals:
\begin{equation}
\gamma(a,z)
=
\displaystyle\int\limits_{0}^{z}{e^{-t} t^{\alpha-1} \, dt},
\quad
\Gamma(a,z)
=
\displaystyle\int\limits_{z}^{\infty}{e^{-t} t^{\alpha-1} \, dt}.
\end{equation}
Parameters $\alpha$ and $z$ are complex numbers and $t^\alpha$ takes its principal value.
Let us remark that the value $\gamma(\alpha,z)$ exists for $\mbox{\rm Re} \, \alpha > 0$
and the value $\Gamma(\alpha,z)$ exists for $|\mbox{\rm arg} \, z| < \pi$. Then,  we have:
$\gamma(a,z) + \Gamma(a,z) = \Gamma(a)$. Analytical continuation can be obtained on
the basis of representation of the $\gamma$ function using series.

\medskip
\noindent
{\cirilsc Mari\ch ev} in \cite{Maricev_78} (page 109, formula 7.116) proved
that for $\mbox{\rm Re} \, z > -1$ the following formula is true:
\begin{equation}
\label{Maricev_78}
\mbox{\rm p.v.}
\displaystyle\int\limits_{0}^{\infty}{e^{-t} \displaystyle\frac{t^z}{t-1} \, dt}
=
\displaystyle\frac{\pi}{e} \mbox{\rm ctg} \, \pi z
-
\Gamma(z) \, \F11(1,1-z,-1),
\end{equation}
where:
\begin{equation}
\F11(a,b,z)
=
\displaystyle\sum\limits_{k=0}^{\infty}{\displaystyle\frac{(a)_{k}}{(b)_{k}}}
\displaystyle\frac{z^{k}}{k!}
\end{equation}
is Kummer's hypergeometric function and where  $(x)_{n} = x (x+1) \ldots
(x+n-1)$, $(x)_{0}=1$ is Pochhammer's symbol. By converting Kummer's hypergeometric
function we get the formula:
\begin{equation}
\label{Konvert_1}
\F11(1,1-z,-1)
=
\displaystyle\frac{(-1)^z}{e}{\big (}\Gamma(1-z)+z\Gamma(-z,-1){\big )}.
\end{equation}
From (\ref{Maricev_78}) and (\ref{Konvert_1}) for $\mbox{\rm Re}\,z\!>\!-1$
the following formula is true:
\begin{equation}
\label{Temme_1}
\mbox{\rm p.v.}
\displaystyle\int\limits_{0}^{\infty}{e^{-t} \displaystyle\frac{t^{z}}{t-1} \, dt}
=
\displaystyle\frac{(-1)^z \Gamma(z+1)\Gamma(-z,-1)}{e}
+
\displaystyle\frac{i \pi}{e}.
\end{equation}
Hence, as one consequence we get representations of functions $K(z)$ and $K_{1}(z)$
via an incomplete gamma function.
\begin{th}
For functions $K(z)$ and $K_{1}(z)$ the following representations are true$:$
\begin{equation}
\label{Repr_K}
K(z)
=
\displaystyle\frac{\mbox{\rm Ei}(1)+i\pi}{e}
+
\displaystyle\frac{(-1)^{z} \Gamma(1+z) \Gamma(-z,-1)}{e}
\end{equation}
and
\begin{equation}
\label{Repr_K1}
K_{1}(z)
=
\displaystyle\frac{(-1)^{z} \pi}{e \sin \pi z}
+
\displaystyle\frac{(-1)^{z} \Gamma(1+z) \Gamma(-z,-1)}{e},
\end{equation}
where the values in the previous formulas, in integer points $z$,
are determined in the sense of Cauchy's principal value.
\end{th}

\noindent
\section{\large \bf \boldmath \hspace*{-7.0 mm}
5. Differential transcendency of functions $K(z)$ and $K_{1}(z)$}

In this section we provide one statement about differential transcendency of some
solutions of functional equations (\ref{K_FE_1}). Namely, using the method of
{\sc Mijajlovi\'c}~\cite{Mijajlovic_98} we can conclude that the following statement is true:
\begin{th}
Let ${\cal M}$ be a differential field of the meromorphic functions $f = f(z)$ over
the connected open set $D \subseteq C \backslash Z^{-}_{0}$. If $g = g(z) \in {\cal M}$
is one solution of a functional equation {\rm (\ref{K_FE_1})}, then $g$ is not a solution
of any algebraic-differential equation over the field of rational functions $C(z)$.
\end{th}
\begin{co}
\label{Posledica_K_K1} Especially the functions $K(z)$ and
$K_{1}(z)$ are not solutions of any algebraic-differential equation
over the field of rational functions $C(z)$.
\end{co}

\bigskip
{\small
\noindent University of Belgrade,
           \hfill (Received September 12, 2003)      \break
\noindent Faculty of Electrical Engineering,   \hfill\break
%           \hfill (Revised ?? ??, ????)             \break
\noindent P.O.Box 35-54, $11120$ Belgrade,     \hfill\break
\noindent Serbia \& Montenegro                 \hfill\break
\noindent {\footnotesize \bf malesevic@kiklop.etf.bg.ac.yu}
\hfill}

\newpage

%%%%%%%%%%%%%%%%%%%%%%%%%%%%%%%%%%%%%%%%%%%%%%%%%%%%%%%%%%%%%%%%%%%%%%%%%%%%%%%%%%%%%%%%%%%%%%%%%

\end{document}